\documentclass[11pt,a4paper]{article}

\usepackage{a4wide,amssymb,graphics,graphicx,textcomp}
\usepackage{enumerate,textcomp,multirow,amsmath}
\usepackage{algorithm}
\usepackage{algorithmic}
\usepackage{url}

\newtheorem{theorem}{Theorem}[section]
\newtheorem{lemma}[theorem]{Lemma}
\newtheorem{remark}[theorem]{Remark}
\newtheorem{definition}[theorem]{Definition}

\newtheorem{example}[theorem]{Example}

\numberwithin{equation}{section}
\numberwithin{table}{section}
\numberwithin{figure}{section}
\newcommand{\proof}{\par\noindent{\bf Proof}. \ignorespaces}
\newcommand{\eproof}{\space
    {\ \vbox{\hrule\hbox{\vrule height1.3ex\hskip0.8ex\vrule}\hrule}}\par}

\newcommand {\mat}  [1] {\left[\begin{array}{#1}}
\newcommand {\rix}      {\end{array}\right]}

\catcode`@=11     
\font\tenex=cmex10 
\newdimen\p@renwd
\setbox0=\hbox{\tenex B} \p@renwd=\wd0 
\def\bmat#1{\begingroup \m@th
  \setbox\z@\vbox{\def\cr{\crcr\noalign{\kern2\p@\global\let\cr\endline}}%
    \ialign{$##$\hfil\kern2\p@\kern\p@renwd&\thinspace\hfil$##$\hfil
      &&\quad\hfil$##$\hfil\crcr
      \omit\strut\hfil\crcr\noalign{\kern-\baselineskip}%
      #1\crcr\omit\strut\cr}}%
  \setbox\tw@\vbox{\unvcopy\z@\global\setbox\@ne\lastbox}%
  \setbox\tw@\hbox{\unhbox\@ne\unskip\global\setbox\@ne\lastbox}%
  \setbox\tw@\hbox{$\kern\wd\@ne\kern-\p@renwd\left[\kern-\wd\@ne
    \global\setbox\@ne\vbox{\box\@ne\kern2\p@}%
    \vcenter{\kern-\ht\@ne\unvbox\z@\kern-\baselineskip}\,\right]$}%
  \null\;\vbox{\kern\ht\@ne\box\tw@}\endgroup}
\catcode`@=12    

\newcommand{\backmatter}{\appendix
\def\chaptermark##1{\markboth{%
\ifnum  \c@secnumdepth > \m@ne  \@chapapp\ \thechapter:  \fi  ##1}{%
\ifnum  \c@secnumdepth > \m@ne  \@chapapp\ \thechapter:  \fi  ##1}}%
\def\sectionmark##1{\relax}}

\makeatletter
\newcommand*{\rom}[1]{\expandafter\@slowromancap\romannumeral #1@}
\makeatother

\def\real{\mathop{\mathrm{Re}}}
\def\imag{\mathop{\mathrm{Im}}}

\newif\ifMDlatex
\catcode`@=11     
\def\MD@us#1{\csname#1style\endcsname}
\def\MD@uf#1{\csname#1font\endcsname}
\def\MD@t{text}
\def\MD@s{script}
\def\MD@ss{scriptscript}
\newdimen\MD@unit
\MD@unit\p@
\def\MD@changestyle#1{
  \relax\MD@unit0.1\fontdimen6\MD@uf{#1}0
  \everymath\expandafter{\the\everymath\MD@us{#1}}
}
\def\MD@dot{$\m@th\ldotp$}
\def\MD@palette#1{\mathchoice{#1\MD@t}{#1\MD@t}{#1\MD@s}{#1\MD@ss}}
\def\MD@ddots#1{{\MD@changestyle{#1}%
  \mkern1mu\raise7\MD@unit\vbox{\kern7\MD@unit\hbox{\MD@dot}}%
  \mkern2mu\raise4\MD@unit\hbox{\MD@dot}%
  \mkern2mu\raise \MD@unit\hbox{\MD@dot}\mkern1mu}}%
\def\MD@iddots#1{{\MD@changestyle{#1}%
  \mkern1mu\raise \MD@unit\hbox{\MD@dot}%
  \mkern2mu\raise4\MD@unit\hbox{\MD@dot}%
  \mkern2mu\raise7\MD@unit\vbox{\kern7\MD@unit\hbox{\MD@dot}}}}%
\def\MD@vdots#1{\vbox{\MD@changestyle{#1}%
    \baselineskip4\MD@unit\lineskiplimit\z@
    \kern6\MD@unit\hbox{\MD@dot}\hbox{\MD@dot}\hbox{\MD@dot}}}%
\ifMDlatex
  \DeclareRobustCommand\ddots{\mathinner{\MD@palette\MD@ddots}}%
  \DeclareRobustCommand\iddots{\mathinner{\MD@palette\MD@iddots}}%
  \DeclareRobustCommand\vdots{\mathinner{\MD@palette\MD@vdots}}%
\else
  \def\ddots{\mathinner{\MD@palette\MD@ddots}}%
  \def\iddots{\mathinner{\MD@palette\MD@iddots}}%
  \def\vdots{\mathinner{\MD@palette\MD@vdots}}%
\fi
\catcode`@=12    

\newcommand {\comment}[1]{} 
\newcommand{\C}{{\mathbb C}}
\newcommand{\R}{{\mathbb R}}

\begin{document}
\title{\hspace{1.3cm} Estimation to structured distances to singularity for matrix pencils with symmetry structures: A linear algebra-based approach}
\author{  Anshul Prajapati $^*$ \qquad Punit Sharma\thanks{A.P. acknowledges the support of the CSIR Ph.D. grant by Ministry of Science \& Technology, Government of India. P.S. acknowledges the support of the DST-Inspire Faculty Award (MI01807-G) by Government of India and FIRP project (FIRP/Proposal Id - 135) by IIT Delhi, India.
 Email: \{maz198078, punit.sharma\}@maths.iitd.ac.in.} \\ 
Department of Mathematics\\
Indian Institute of Technology Delhi\\
Hauz Khas, 110016, India
}

\date{}

\maketitle

\begin{abstract}
We study the structured distance to singularity for a given regular matrix pencil $A+sE$, where $(A,E)\in \mathbb S \subseteq (\C^{n,n})^2$. This includes
Hermitian, skew-Hermitian, $*$-even, $*$-odd, 
$*$-palindromic, T-palindromic, and dissipative Hamiltonian pencils. We present a
purely linear algebra-based approach to 
derive explicit computable formulas for the distance to the nearest structured pencil $(A-\Delta_A)+s(E-\Delta_E)$ such that $A-\Delta_A$ and $E-\Delta_E$ have a common null vector. We then obtain a family of computable lower bounds for the unstructured and structured distances to singularity. 
Numerical experiments suggest that in many cases, there is a significant difference between structured and unstructured distances.
 This approach extends to structured matrix polynomials with higher degrees.
\end{abstract}

{\textbf{ keyword}}
structured matrix pencils, eigenvalue backward errors, distance to singularity, structured distance to singularity, distance to common null space, differential-algebraic system, dissipative Hamiltonian system


\section{Introduction}
Let $L(s) = A+sE$ be a matrix pencil where $A$ and $E$ are $n \times n$ complex matrices,
then $L(s)$ is {\em regular} if $A+\lambda E$ is invertible for some $\lambda \in \C$. A matrix pencil that is
not regular is often called a {\em singular matrix pencil}. The matrix pencil $L(s)$ is said to be
\emph{structured} if $(A,E)$ belongs to a special subset $\mathbb S$ of $(\C^{n,n})^2$.
We consider Hermitian, skew-Hermitian, $*$-even, $*$-odd, 
$*$-palindromic,  T-palindromic, and dissipative Hamiltonian structures~\cite{Lan66a,MehMW20,Meh91, Sch06,SchM95}.

Given a regular structured matrix pencil $A+s E$, where $(A,E)\in \mathbb S$, the problem of
determining the nearest singular pencil $(A-\Delta_A)+s (E-\Delta_E)$, where $(\Delta_A,\Delta_E)\in \mathbb S$
with respect to some specific norm is called the {\em structured distance to
singularity} for matrix pencil $A+s E$. If $\mathbb S=(\C^{n,n})^2$, then it is called the
{\em unstructured distance to singularity} (a.k.a. {\em distance to singularity
}) for pencil $A+s E$.
Due to its great importance in many engineering applications, this problem has
attracted much work and is still an open problem \cite{BerGTWW19,ByeHM98,GugLM17,MehMW15,MehMW17}.
As an example, in the linear time-invariant descriptor systems
\begin{equation}\label{eq:daesys}
E\dot{x(t)}=Ax(t)+Bu(t), \quad y(t)=Cx(t)+Du(t),
\end{equation}
if the matrix pencil $A+s E$ is singular, then the initial value problem
of solving~\eqref{eq:daesys} with a consistent initial value $x_0$ is not
solvable or the solution is not unique. Thus the system equation~\eqref{eq:daesys} is not
well-posed if $A+sE$ is singular. Also, an unfortunate choice of linear feedback control $u=Fx$ or $u=Fy$
may make
the system singular, or nearly singular.
Further, eigenvalues and eigenvectors of a regular matrix pencil may vary discontinuously in the
neighbourhood of a singular pencil. Thus, in such cases, the sensitivity of the problem is closely
related to the nearest ill-posed problem \cite{Dem87}, and some special numerical methods are required \cite{DemK88}.
If the matrix pencil with additional structure is considered, then the use of structure-preserving algorithms is advisable.

Although a complete characterization of the distance to singularity is still open, many heuristic bounds are available in the literature \cite{BerGTWW19,ByeHM98,GugLM17,MehMW15,MehMW17}. The situation can be different if $E$ and $A$ follow some symmetry structure and one considers the structured distance to singularity~\cite{MehMW17,MehMW20}. Recently~\cite{MehMW20}, considered the structured distance to singularity for dissipative Hamiltonian (DH) pencils, $sE+(J-R)Q$, where $J^T=-J$, and $E^TQ$ and $R$ are positive semidefinite. The authors in~\cite{MehMW20} characterize the structured distance to singularity for DH pencils in terms of null spaces of some structured matrices for the case when $Q=I_n$. Some special types of DH pencils fall into subclasses of even pencils and Hermitian pencils. For example, DH pencils with $R=0$ and $Q=I_n$, is an odd pencil, and  similarly with $J=0$ and $Q=I_n$, is a Hermitian pencil. 

 Another motivation to study the structured distance to singularity comes from the passivity analysis of the system~\eqref{eq:daesys}. Such a system is passive if the even pencil
 \[
 L(s)=sM+N:=s\mat{ccc}0 & E & 0\\-E^* &0 &0 \\ 0&0&0\rix 
 + \mat{ccc}0&A&B\\A^* & 0& C^* \\ B^* & C& D+D^*\rix
 \]
is regular, has no purely imaginary eigenvalues, and is of the index at most one~\cite{FreJ04}. Thus, it is important to know for which perturbation system becomes singular. This singularity can be associated with
a structure-preserving perturbation
\[
s\Delta_M+\Delta_N:=s\mat{ccc}0 & \Delta_E & 0\\-\Delta_E^* &0 &0 \\ 0&0&0\rix 
 + \mat{ccc}0&\Delta_A&\Delta_B\\\Delta_A^* & 0& \Delta_C^* \\ \Delta_B^* & \Delta_C& \Delta_D+\Delta_D^*\rix
\] 
to the pencil $sM+N$ that preserves the even structure and the zero block structure of the pencil. So if this perturbation is smaller than the smallest structured perturbation that makes the system singular, then the system remains regular. Thus structured distance to singularity gives an upper bound to the passivity radius. This is another motivation to study a structured distance. A much harder structure that we do not address is to preserve the even and the zero blocks structure of pencil $sM+N$ which would be more appropriate for the passivity of the system~\eqref{eq:daesys}~\cite{MehV21}. 

In the simplest case, a pencil is singular due to a common left or right null vector. Thus the {\em distance to a common null space}, i.e., the distance to the nearest pencil 
$(A-\Delta_A)+s(E-\Delta_E)$ such that
$(A-\Delta_A)$ and $(E-\Delta_E)$ have a common null vector, gives an upper bound to the distance to singularity. If the pencil is structured, then the {\em structured distance to a common null space} is more appropriate. 

\subsection{Contribution and outline of the paper}

In Section~\ref{sec:prel}, we introduce various structured distances and recall some preliminary results from literature that will be useful for the main results of the paper. 

In Section~\ref{sec:dist_null}, we present a purely
linear algebra-based approach to obtain computable formulas for the structured distance to a common null space for Hermitian, skew-Hermitian, $*$-even, $*$-odd, $*$-palindromic, T-palindromic, and DH pencils. The formulas have been obtained with respect to the measure $\sqrt{{\|A\|}^2+{\|E\|}^2}$, where ${\|\cdot\|}$ is the matrix spectral norm. For DH pencils of the form $sE+(J-R)$, the structured distance to a common null space is obtained with respect to structure-preserving skew-Hermitian perturbations to $J$ and negative semidefinite perturbations to $R$ and $E$. This distance
is shown to be equal to the structured distance to singularity when the perturbations in $R$ and $E$ are restricted to be of rank one. We note that this distance is different from the one considered in~\cite{MehMW20}, where the DH structure was preserved under symmetric and possibly indefinite perturbations to $R$ and $E$. 

The techniques for deriving the structured distance to a common null space are further extended to the matrix polynomials of degree more than one. These results are briefly discussed in Section~\ref{sec:highdeg}.

An equivalent characterization for singular pencils is that a pencil $A+sE,~A, E\in \C^{n,n}$ is singular if and only if there exists $(\lambda_1,\ldots,\lambda_{n+1}) \in \C^{n+1}$ such that
$\lambda_i \neq \lambda_j$ if $i \neq j$, and $\operatorname{det}(A+\lambda_i E)=0$ for $i=1,\ldots,n+1$.
Thus the distance to singularity problem can be reformulated in terms of finding the nearest pencil that has some pre-specified distinct $n+1$ points in the complex plane as eigenvalues.
Motivated by these considerations, we obtain a family of lower bounds for the unstructured and structured distances to singularity in Section~\ref{sec:dist_sing}.
These bounds involve minimizing the largest eigenvalue of some parameter-dependent Hermitian matrix which can be computed by using a
suitable optimization technique.  

In Section~\ref{sec:numresults}, we present some numerical experiments to highlight the significance of structure-preserving and arbitrary perturbations on the various distances under consideration.

\section{Notation and preliminaries}\label{sec:prel}

In the following, we denote the spectral norm of a matrix or a vector by $\|\cdot\|$, the Frobenius norm of a matrix by ${\|\cdot\|}_F$, the smallest and the largest eigenvalues of a Hermitian matrix $A$ respectively by $\lambda_{\min}(A)$ and $\lambda_{\max}(A)$, the Moore-Penrose pseudoinverse of a matrix or a vector $X$ by $X^\dagger$ and the smallest singular value of a matrix $A$ by $\sigma_{\min}(A)$. We use
$\text{Herm}(n)$ and $\text{SHerm}(n)$ respectively to denote the set of Hermitian  and skew-Hermitian matrices of size $n \times n$ and
$\mathcal{E}_{\lambda}(A)$ to denote the eigenspace of the matrix $A$ corresponding to the eigenvalue $\lambda$. The second-largest eigenvalue of a Hermitian matrix $A$ is denoted by $\lambda_2(A)$ and the second-largest singular value of a matrix $A$ is denoted by $\sigma_2(A)$. We use the notation $A\succeq 0$ ($A \preceq 0$) and $A\succ 0$ ($A\prec 0$) if the matrix is Hermitian positive semidefinite (negative semidefinite) and Hermitian positive definite (negative definite), respectively. By $S^c$ we denote the set complement of the set $S$ within a larger set. 

\begin{definition}
Consider a regular matrix pencil $L(s)=A+s E$, where $(A, E) \in \mathbb{S} \subseteq (\C^{n, n})^2$.
\begin{enumerate}
\item The \emph{structured distance to a common null space} with respect to structured perturbations to both $A$ and $E$ is defined by
\begin{equation}\label{nullAE}
    \delta_0^{\mathbb S}(A,E):= \inf\big\{\sqrt{\|\Delta_A\|^2+\|\Delta_E\|^2} \; : \;(\Delta_A, \Delta_E) \in \mathbb{S}, \; \text{\rm ker}(A-\Delta_A)\cap \text{\rm ker}(E-\Delta_E) \neq \{0\}\big\}.
\end{equation}
\begin{itemize}
    \item When $\mathbb{S} = (\C^{n, n})^2$, then $\delta_0^{\mathbb S}(A,E)$ is called the {\em unstructured distance to a common null space} and  we denote it by
    $\delta_0(A,E) := \delta_0^{\mathbb S}(A,E)$.
    \item $\delta_0^{\mathbb S}(A,E) \leq \sqrt{\|A\|^2+\|E\|^2}$ because ${\rm ker}(A-\Delta_A) \cap {\rm ker}(E - \Delta_E) \neq \{0\}$ holds trivially for $\Delta_A = A$ and $\Delta_E = E.$
\end{itemize}
\item The \emph{structured distances to a common null space} with respect to structured perturbations to only one of $A$ and $E$ are respectively defined by
\begin{equation}\label{nullAA}
    \delta_0^{\mathbb S}(A):= \inf\big\{\|\Delta_A\| \; : \;\Delta_A\in \C^{n,n},\;(A-\Delta_A, E) \in \mathbb{S}, \; \text{\rm ker}(A-\Delta_A)\cap \text{\rm ker}(E) \neq \{0\}\big\};
\end{equation}
and
\begin{equation}\label{nullEE}
    \delta_0^{\mathbb S}(E):= \inf\big\{\|\Delta_E\| \; : \;\Delta_E\in \C^{n,n},\;(A,E-\Delta_E) \in \mathbb{S}, \; \text{\rm ker}(A)\cap \text{\rm ker}(E-\Delta_E) \neq \{0\}\big\}.
\end{equation}
The analogous unstructured distances are denoted by $ \delta_0(A)$ and $\delta_0(E)$.
\item The \emph{structured distance to singularity}, $\delta^{\mathbb S}(A,E)$ with respect to structured perturbations to both $A$ and $E$ is
the smallest perturbation $(\Delta_A, \Delta_E) \in \mathbb S$ such that $(A-\Delta_A)+z(E-\Delta_E)$ is singular. More precisely,
\begin{equation}\label{sdistance}
    \delta^{\mathbb S}(A,E):= \inf\{\sqrt{\|\Delta_A\|^2+\|\Delta_E\|^2} \; : \;(\Delta_A, \Delta_E) \in \mathbb S, \; \det((A-\Delta_A)+\lambda(E-\Delta_E)) = 0 \; {\rm for\,all} \; \lambda \in \C\}.
\end{equation}
\begin{itemize}
    \item When $\mathbb{S} = (\C^{n, n})^2$, then $\delta^{\mathbb S}(A,E)$ is called the {\em unstructured distance to singularity} and  we denote it by
    $\delta(A,E) := \delta^{\mathbb S}(A,E)$.
    \item $\delta^{\mathbb S}(A,E) \leq \sqrt{\|A\|^2+\|E\|^2}$ because for $(\Delta_A, \Delta_E) = (A, E)$ the determinant condition, $ \det((A-\Delta_A)+\lambda(E-\Delta_E)) = 0 $, is satisfied for all $\lambda$.
\end{itemize}
\item The \emph{structured distances to singularity} with respect to structured perturbations to only one of $A$ and $E$ are respectively defined by
\begin{equation}\label{nullA}
    \delta^{\mathbb S}(A):= \inf\big\{\|\Delta_A\| \; : \;\Delta_A\in \C^{n,n},\;(A-\Delta_A, E) \in \mathbb{S}, \; \det((A-\Delta_A)+\lambda E) = 0\; {\rm for\, all} \; \lambda \in \C\big\};
\end{equation}
and
\begin{equation}\label{nullE}
    \delta^{\mathbb S}(E):= \inf\big\{\|\Delta_E\| \; : \;\Delta_E\in \C^{n,n},\;(A,E-\Delta_E) \in \mathbb{S}, \; \det(A+\lambda(E-\Delta_E)) = 0\; {\rm for\,all} \; \lambda \in \C\big\}.
\end{equation}
The analogous unstructured distances are denoted by $ \delta(A)$ and $\delta(E)$. Clearly, $\delta(E)=\infty$ when $A$ is invertible.
\end{enumerate}
 If the perturbations are restricted to be of rank one, then we denote the above distances by adding an index 1, i.e., we write $\delta_1$ for the corresponding distance. For example, $\delta_{10}^\mathbb S(A,E)$ 
denotes the structured distance to a common null space with respect to
structure preserving rank one perturbations to $A$ and $E$. 
\end{definition}

In the following, we recall some mapping results from the literature which will be crucial for us in deriving various structured distances.

\begin{lemma}{\rm\cite{MacMT08}}\label{l1}
Let $x \in \mathbb{K}^n\setminus \{0\}$ and $b \in \mathbb{K}^n$ where $\mathbb{K} \in \{\R, \C\}$. Then there exists $ \Delta \in \mathbb{K}^{n\times n}$ such that  $\Delta x=b$ iff $\Delta$ is of the form
\[
\Delta = bx^{\dagger} + Z(I_n - xx^{\dagger}),
\]
where $Z \in \mathbb{K}^{n\times n}$ is arbitrary. Furthermore, the minimal spectral norm of such $\Delta$ is $\frac{\|b\|}{\|x\|}$ and is attained by $\Delta = bx^{\dagger}.$
\end{lemma}

\begin{lemma}{\rm\cite{MacMT08}}\label{l2}
Let $x, y \in \C^n,\; x \neq 0$. Then there exists a matrix $H \in {\rm Herm}(n)$ such that $Hx = y$ if and only if $\imag{(x^*y)} = 0$. If the latter condition is satisfied, then
\[
\min\big\{\|H\|\; :\; H \in {\rm Herm}(n), \; Hx=y\big\} = \frac{\|y\|}{\|x\|}.
\]
\end{lemma}

\begin{lemma}{\rm\cite[Theorem 2.5]{BorKMS15}}\label{l3}
Let $x,y,z\in\C^n$ with $x\neq 0$, and $\star \in \{*,T\}$. Then there exists a matrix $\Delta \in \C^{n\times n}$
such that $\Delta x=y$ and $\Delta ^\star x=z$ if and only if
$x^\star y=z^\star x$. If the latter condition is satisfied, then
\begin{equation}\label{eq:minpalin}
\min\Big\{ \|\Delta \|\,:\, \Delta \in\C^{n\times n}, \,
  \Delta x=y,\Delta ^\star x=z\Big\}=
\max\left\{\frac{\|y\|}{\|x\|},\frac{\|z\|}{\|x\|}\right\}.
\end{equation}
\end{lemma}

\begin{lemma}{\rm \cite[Theorem 2.3]{MehMS16}}\label{lemdh2}
Let $x,y \in \C^{n}\setminus \{0\}$. Then there exists a matrix $\Delta \in \C^{n,n}$
such that $\Delta \preceq 0$ and $\Delta x=y$ if and only if $x^*y < 0$. 
If the latter condition is satisfied then
\[
\min\left\{\|\Delta\|:~\Delta \preceq 0,\, \Delta x=y\right\}=\frac{{\|y\|}^2}{|x^*y|}
\]
and the minimum is attained by the rank one matrix $\tilde \Delta=\frac{1}{x^*y}yy^*$.
\end{lemma}
\begin{lemma}{\rm \cite[Lemma 4.1]{MehMS16}}\label{lemdh4}
Let $H \in \C^{n,n}$ be such that $H^*=H \succeq 0$. Suppose $x \in \C^n$
such that $Hx \neq 0$ and set $\Delta_H=-\frac{(Hx)(Hx)^*}{x^*Hx}$,
then $H+\Delta_H$ is Hermitian positive semidefinite. 
\end{lemma}
The following lemma will be useful in deriving explicit formula for the structured distance to a common null space for palindromic pencils.
\begin{lemma}{\rm \cite[Proposition 6.1]{Kar11}}\label{l4}
Let $H_0, H_1 \in$ {\rm Herm}($n$) and $H(t)=H_0 + tH_1,\, t \in \R$. Suppose the function $t\mapsto \lambda_{\min}(H(t)),\, t \in \R$ attains a local extremum at $t_0$.
 Then there exists a vector $v \in \mathcal{E}_{\lambda_{\min}}(H(t_0))$ such that  $\|v\|=1$ and $v^*H_1v = 0$.
\end{lemma}

A real version of the following lemma was proved in~\cite{MehMW20}. We state the result for a complex DH pencil. 
\begin{lemma}\label{lemdh1}
Let $L(s)=sE-(J-R)$ be a DH pencil, where $J,R,Q \in \C^{n,n}$ such that $J^*=-J$, $R \succeq 0$, and $E \succeq 0$. Then 
\[
\text{\rm ker}(E)\cap \text{\rm ker}(J-R)=\text{\rm ker}(E) \cap \text{\rm ker}(R) \cap \text{\rm ker}(J). 
\]
\end{lemma}

\begin{lemma}{\rm \cite[Theorem 5]{MehMW20}}\label{lemdh3}
Let $L(s)=sE+(J-R)$ be a real DH pencil, where $J,R,Q \in \R^{n,n}$ such that 
$R\succeq 0$, $E \succeq 0$, and $J^T=-J$. Then $L(s)$ is singular if and only if 
\[
\text{\rm ker}(E) \cap \text{\rm ker}(R) \cap \text{\rm ker}(J) \neq \{0\}.
\]
\end{lemma}

\section{Structured distances to a common null space}\label{sec:dist_null}

In this section, we derive an explicit computable formula for the structured distance to a common null space $\delta_0^{\mathbb S}(A,E)$ defined in~\eqref{nullAE} for the structures under consideration.

\subsection{Unstructured pencils}\label{subsec:unstrnullspace}

The following result gives an explicit formula for the unstructured distance 
$\delta_0(A,E)$.

\begin{theorem}\label{thm:unstrnull}
Let $L(s)=A+s E$ be a regular matrix pencil , where $(A, E) \in (\C^{n, n})^2$. Then
\[
\delta_0(A,E)= \sqrt{\lambda_{\min}(A^*A+E^*E)}.
\]
\end{theorem}
\proof First notice that for any $(\Delta_A,\Delta_E)\in (\C^{n,n})^2$, $\text{\rm ker}(A-\Delta_A) \cap \text{\rm ker}(E-\Delta_E) \neq \{0\}$ if and only if
there exists $ v \in \C^{n}\setminus \{0\}$ such that $ (A-\Delta_A)v=0$ and $(E-\Delta_E)v=0$ if and only if $\Delta_A v = Av$  and $\Delta_E v = Ev$.
Using this in the definition~\eqref{nullAE}, we get
\begin{equation}\label{uneq1}
\delta_0(A,E)^2 = \inf\left \{\|\Delta_A\|^2+\|\Delta_E\|^2 \; : \;\Delta_A, \Delta_E \in \C^{n \times n},\,v\in \C^{n}\setminus\{0\},\,\Delta_Av=Av, \Delta_E v=Ev \right\}.
\end{equation}
Thus in view of Lemma~\ref{l1}, a minimal norm $(\Delta_A,\Delta_E)$ for any fixed nonzero vector $v \in \C^n$ such that $\Delta_Av=Av$ and $\Delta_E v=Ev$, is given by
\[
\|\Delta_A\|^2 +\|\Delta_E\|^2 = \frac{\|Av\|^2}{\|v\|^2}+\frac{\|Ev\|^2}{\|v\|^2} = \frac{v^*(A^*A+E^*E)v}{v^*v},
\]
and is attained by rank-1 perturbations $\Delta_A=Avv^\dagger$ and
$\Delta_E=Evv^\dagger$.
This yields from~\eqref{uneq1} that
\begin{eqnarray*}
        \delta_0(A,E)^2
        & =& \inf \left\{\frac{v^*(A^*A + E^*E)v}{v^*v} \;: \; v \in \C^{n}\setminus \{0\}\right\}\\
        & =& \lambda_{\min}(A^*A + E^*E).
\end{eqnarray*}
\eproof
The distance to a common null space for $A+sE$ with respect to perturbations only in the matrix $A$ is obtained in the following result.
\begin{theorem}\label{ustrnullA}
Let $L(s)=A+sE$ be a regular matrix pencil, where $A,E\in \C^{n,n}$. Then
\[
\delta_0(A)=
\begin{cases}
\infty &\text{if \,} E \,\text{is\, invertible},\\
\sigma_{\min}(AU) &\text{otherwise },
\end{cases}
\]
where the columns of $U$ form an orthonormal basis for ${\rm ker}(E)$.
\end{theorem}
\proof From~\eqref{nullAA}, if $E$ is invertible, then clearly $\delta_0(A)=\infty$. Thus assume that $\text{\rm ker}(E)\neq \{0\}$ and let $k=\text{dim}(\text{\rm ker}(E))$. Then
\begin{eqnarray}\label{eq:nillA1}
    \delta_0(A)^2&=& \inf\left\{\|\Delta_A\|^2 \; : \;\Delta_A\in \C^{n,n},\;\text{\rm ker}(A-\Delta_A)\cap \text{\rm ker}(E) \neq \{0\}\right\}\nonumber\\
    &=& \inf\left\{\|\Delta_A\|^2 \; : \;\Delta_A\in \C^{n,n},\; x\in \C^n\setminus \{0\},\; (A-\Delta_A)x=0,\;Ex=0\right\}\nonumber\\
    &=& \inf\left\{\frac{\|Ax\|^2}{\|x\|^2} \; : \; x\in \C^n\setminus \{0\},\; Ex=0\right\},
\end{eqnarray}
where the last equality follows by using the minimal norm mapping from Lemma~\ref{l1}. Let $U \in \C^{n,k}$ be such that its columns form an orthonormal basis for ker$(E)$. Then for any
nonzero $x\in \text{\rm ker}(E)$, there exists $\alpha \in \C^{k}\setminus \{0\}$ such that $x=U\alpha$. Using this in~\eqref{eq:nillA1}, we have
\[
    \delta_0(A)^2= \inf\left\{\frac{\|AU\alpha\|^2}{\|\alpha\|^2} \; : \; \alpha\in \C^k\setminus \{0\}\right\}= (\sigma_{\min}(AU))^2,
\]
because $\|U\alpha\|=\|\alpha\|$ for any $\alpha\in \C^k\setminus \{0\}$ as $\|\cdot\|$ is unitarily invariant.
\eproof

Analogously, we can obtain $\delta_0(E)$ as following:
\[
\delta_0(E)=
\begin{cases}
\infty &\text{if \,} A \,\text{is\, invertible},\\
\sigma_{\min}(EV) &\text{otherwise },
\end{cases}
\]
where the columns of $V$ form an orthonormal basis for $\text{\rm ker}(A)$.

\begin{remark}{\rm
Note that in Theorem~\ref{thm:unstrnull} the distance  $\delta_0(A,E)$ is attained by rank one perturbations $\Delta_A=Avv^\dagger$ and $\Delta_E=Evv^\dagger$, where $v$ is an eigenvector corresponding to the eigenvalue $\lambda_{\min}(A^*A+E^*E)$ of $A^*A+E^*E$.
Therefore, we have that
\[
\delta_0(A,E)={\delta_1}_0(A,E).
\]
A perturbation that attains the distance $\delta_0(A)$ in Theorem~\ref{ustrnullA} can be constructed as follows: take $u$ to be a right singular vector corresponding to the singular value $\sigma_{\min}(AU)$ of matrix $AU$, set $x=Uu$ and consider the perturbation $\Delta_A=Axx^\dagger$. Then this $\Delta_A$ is of rank one such that $(A-\Delta_A)$ and $E$ has a common null vector, and $\|\Delta_A\|=\sigma_{\min}(AU)$. This also implies that
\[
\delta_0(A)={\delta_1}_0(A).
\]
}
\end{remark}

\subsection{Hermitian and related structures}\label{subsec:null_herm}

In this section, we consider Hermitian pencils $L(s)$, i.e., $L(s)=A+sE$ where $(A,E) \in \mathbb S=(\text{Herm}(n))^2$. We first show that
for Hermitian pencils the structured distance to common null space is equal to the unstructured distance.
\begin{theorem}\label{thm:hermnull}
Let $\mathbb S=(\text{\rm Herm}(n))^2$ and let  $L(s)=A+s E$ be a regular Hermitian pencil , where $(A, E) \in \mathbb S$. Then
\[
\delta_0^{\text{\rm Herm}}(A,E)= \sqrt{\lambda_{\min}(A^*A+E^*E)}=\delta_0(A,E).
\]
\end{theorem}
\proof By definition~\eqref{nullAE}, we have
\begin{eqnarray}\label{errorh}
\delta_0^{\text{Herm}}(A,E)^2 =& \inf\Big\{\|\Delta_A\|^2 + \|\Delta_E\|^2\; :\; \Delta_A, \Delta_E \in \text{Herm}(n),\, \text{\rm ker}(A-\Delta_A) \cap
\text{\rm ker}(E-\Delta_E) \neq \{0\}\Big\}.\nonumber\\
=&\hspace{-3.9cm}\inf\Big\{\|\Delta_A\|^2 + \|\Delta_E\|^2\; :\; \Delta_A, \Delta_E \in \text{Herm}(n),\; v \in \C^n\setminus\{0\},\nonumber\\
&\hspace{-2cm} (A-\Delta_A)v=0, \; (E-\Delta_E)v=0\Big\} \nonumber\\
=&\inf\Big\{\|\Delta_A\|^2 + \|\Delta_E\|^2\; :\; \Delta_A, \Delta_E \in \text{Herm}(n),\; v \in \C^n\setminus\{0\},\,\Delta_Av=Av, \; \Delta_Ev=Ev\Big\},\nonumber \\
\end{eqnarray}
where the last equality follows by similar arguments as used in~\eqref{uneq1}.
In view of Lemma~\ref{l2}, a minimal norm $(\Delta_A,\Delta_E) \in (\text{Herm}(n))^2$ for a fixed nonzero vector $v \in \C^n$ such that $\Delta_Av=Av$ and $\Delta_E v=Ev$, is given by
\[
\|\Delta_A\|^2 +\|\Delta_E\|^2 = \frac{\|Av\|^2}{\|v\|^2}+\frac{\|Ev\|^2}{\|v\|^2} = \frac{v^*(A^*A+E^*E)v}{v^*v}.
\]
This quantity is minimal in norm for a fixed vector $v$ and now we have to minimize this over all possible vectors $v \in \C^n$ for which there
exists $(\Delta_A,\Delta_E) \in (\text{Herm}(n))^2$ such that $\Delta_Av=Av$ and $\Delta_E v=Ev$. By Lemma~\ref{l2}, this happens only for
vectors $v$ satisfying $\imag{v^*Av}=0$ and $\imag{v^*Ev}=0$. Since $A, E \in \text{Herm}(n)$, the conditions
$\imag{v^*Av}=0$ and $\imag{v^*Ev}=0$  trivially hold for every vector $v$. Using this in~\eqref{errorh} we obtain
\begin{eqnarray*}
        \delta_0^{\text{Herm}}(A,E)^2
        & =& \inf \left\{\frac{v^*(A^*A + E^*E)v}{v^*v} \;: \; v \in \C^{n}\setminus \{0\}\right\}\\
        & =& \lambda_{\min}(A^*A + E^*E)\\
        & =&  \delta_0(A,E)^2,
\end{eqnarray*}
where the last equality is due to Theorem~\ref{thm:unstrnull}.
\eproof

\begin{remark}\label{rem:nullherm}{\rm
In view of Theorems~\ref{ustrnullA} and~\ref{thm:hermnull}, an analogous result can be obtained for the structured distances $\delta_0^{\text{Herm}}(A)$ and $\delta_0^{\text{Herm}}(E)$, i.e., we have
 \[
\delta_0^{\text{Herm}}(A)=\delta_0(A) \quad \text{and}\quad \delta_0^{\text{Herm}}(E)=\delta_0(E).
\]
}
\end{remark}

Note that some other cases of structured pencils can be converted to the Hermitian case and thus using Theorem~\ref{thm:hermnull}, we can derive
the distance to a common null space for these structures also. Indeed, if $L(s)$ is a skew-Hermitian pencil, then $P(s) = iL(s)$ is a Hermitian pencil. Similarly, if $L(s)$ is $*$-even
or $*$-odd, then one may instead consider the Hermitian pencils $P(s) = L(is)$ or $P(s) = iL(is)$. Let us denote $\mathbb S_e$ as the
set of  $*$-even matrix pairs $(\Delta_0,\Delta_1)$ such that $\Delta_0 \in \text{Herm}(n),\;\Delta_1 \in \text{SHerm}(n)$, and $\mathbb S_o$ as the set of  $*$-odd matrix pairs $(\Delta_0,\Delta_1)$ such that $\Delta_0 \in \text{SHerm}(n),\;\Delta_1 \in \text{Herm}(n)$. Then we have the following theorem.

\begin{theorem}\label{thm:shermnull}
Let $\mathbb S\in \{({\rm SHerm}(n))^2,\mathbb S_e,\mathbb S_o\}$ and let  $L(s)=A+s E$ be a regular structured  pencil, where $(A, E) \in \mathbb S$. Then
\[
\delta_0^{\text{\rm SHerm}}(A,E)=\delta_0^{\text{\rm Herm}}(iA,iE)= \sqrt{\lambda_{\min}(A^*A+E^*E)}=\delta_0(A,E),
\]
\[
\delta_0^{\mathbb S_e}(A,E)=\delta_0^{\text{\rm Herm}}(A,iE)= \sqrt{\lambda_{\min}(A^*A+E^*E)}=\delta_0(A,E),
\]
and
\[
\delta_0^{\mathbb S_o}(A,E)=\delta_0^{\text{\rm Herm}}(iA,-E)= \sqrt{\lambda_{\min}(A^*A+E^*E)}=\delta_0(A,E).
\]
\end{theorem}

%
%

We note that Remark~\ref{rem:nullherm} can also be generalized for skew-Hermitian, $*$-even, and $*$-odd pencils.

\subsection{Palindromic structures}\label{sec:pal_null}

Following the terminology in~\cite{BorKMS15}, we use the term $\star$-palindromic where $\star=*$ or $\star=T$, whenever a statement
is valid for both $*$-palindromic and T-palindromic structures. In the following theorem, unlike the Hermitian structure,
we show that the structured distance to a common null space for $\star$-palindromic pencil is different from the unstructured one.
\begin{theorem}\label{thm:res_null_pal}
Let $L(s)=A+sA^\star$ be a regular $\star$-palindromic pencil, where $A\in \C^{n,n}$. Then
\begin{equation}
\delta_0^{\text{\rm pal}_*}(A,A^*)^2 = 2\sup_{\gamma \in [0,1]}\lambda_{\min}\left(AA^* + \gamma(A^*A - AA^*)\right),
\end{equation}
when $\star=*$, and
\begin{equation}
\delta_0^{\text{\rm pal}_T}(A,A^T)^2 = 2\sup_{\gamma \in [0,1]}\lambda_{\min}\left(\bar{A}A^T + \gamma(A^*A - \bar{A}A^T)\right),
\end{equation}
when $\star=T$.
\end{theorem}
\proof
By definition~\eqref{nullAE}, we have
\begin{eqnarray}\label{errorpal}
\delta_0^{\text{pal}_\star}(A,A^\star) &=& \inf\Big\{\sqrt{\|\Delta_A\|^2 + \|\Delta_A^\star\|^2}\; :\; \Delta_A\in \C^{n,n},\, \text{\rm ker}(A-\Delta_A) \cap
\text{\rm ker}(A^\star-\Delta_A^\star) \neq \{0\}\Big\}.\nonumber\\
&=&\sqrt{2}\inf\Big\{\|\Delta_A\|\; :\; \Delta_A\in \C^{n,n},\; v \in \C^n\setminus\{0\},\;(A-\Delta_A)v=0, \; (A^\star-\Delta_A^\star)v=0\Big\} \nonumber\\
&=&\sqrt{2}\inf\Big\{\|\Delta_A\|\; :\; \Delta_A\in \C^{n,n},\; v \in \C^n\setminus\{0\},\,\Delta_Av=Av, \; \Delta_A^\star v=A^\star v\Big\},
\end{eqnarray}
since $\|\Delta\|=\|\Delta^\star\|$ for any $\Delta \in \C^{n,n}$. In view of Lemma~\ref{l3}, for any fixed $v\in \C^n\setminus\{0\}$, there exists
a $\Delta_A\in \C^{n,n}$ such that $\Delta_A v=Av$ and $\Delta_A^\star v=A^\star v$ if and only if $v^\star A v=(A^\star v)^\star v$. Also, such
a $\Delta_A$ which is also minimal with respect to spectral norm satisfies
\[
\|\Delta_A\|=\max\left\{\frac{\|Av\|}{\|v\|},\frac{\|A^\star v\|}{\|v\|}\right\}.
\]
Note that the necessary and sufficient condition $v^\star A v=(A^\star v)^\star v$ holds trivially for any $v\in \C^n$. Thus~\eqref{errorpal} yields
\begin{equation}\label{p1}
\delta_0^{\text{pal}_\star}(A,A^\star)^2 = 2\inf\left\{\max\left\{\frac{\|Av\|^2}{\|v\|^2},\frac{\|A^\star v\|^2}{\|v\|^2}\right\} \; : \; v\in \C^n\setminus \{0\}\right\}.
\end{equation}
Observe that for any $v\in \C^n\setminus \{0\}$ and $\gamma \in [0,1]$, we have
\begin{equation}\label{p2}
    \max\left\{\frac{\|Av\|^2}{\|v\|^2},\frac{\|A^\star v\|^2}{\|v\|^2}\right\} \geq \gamma \frac{\|Av\|^2}{\|v\|^2} + (1-\gamma)\frac{\|A^\star v\|^2}{\|v\|^2}.
\end{equation}
After taking the infimum in~\eqref{p2} over all $v\in \C^n\setminus \{0\}$, we get
\[
\delta_0^{\text{pal}_\star}(A,A^\star)^2 \geq 2 \inf\left\{ \gamma \frac{\|Av\|^2}{\|v\|^2} + (1-\gamma)\frac{\|A^\star v\|^2}{\|v\|^2} \;:\; v\in \C^n\setminus \{0\}\right\},
\]
which is true for every $\gamma \in [0,1]$. This implies that
\begin{eqnarray}
\delta_0^{\text{pal}_\star}(A,A^\star)^2 &\geq& 2 \sup_{\gamma \in [0,1]}\inf\left\{ \gamma \frac{\|Av\|^2}{\|v\|^2} + (1-\gamma)\frac{\|A^\star v\|^2}{\|v\|^2} \;:\; v\in \C^n\setminus \{0\}\right\} \label{temp1}\\
&=& 2\sup_{\gamma \in [0,1]}\inf\left\{ \gamma \frac{v^*A^*Av}{v^*v} + (1-\gamma)\frac{v^*(A^\star)^*A^\star v }{v^*v} \;:\; v\in \C^n\setminus \{0\}\right\}  \nonumber \\
&=& 2\sup_{\gamma \in [0,1]}\inf\left\{  \frac{v^*\left((A^\star)^*A^\star+\gamma (A^*A-(A^\star)^*A^\star)\right)v}{v^*v} \;:\; v\in \C^n\setminus \{0\}\right\}  \nonumber \\
&=& 2\sup_{\gamma \in [0,1]} \lambda_{\min}\left((A^\star)^*A^\star+\gamma (A^*A-(A^\star)^*A^\star)\right).\label{temp2}
\end{eqnarray}
The proof is complete if we show that equality holds in~\eqref{temp1} and thus a computable formula for $\delta_0^{\text{pal}_\star}(A,A^\star)$ is given by~\eqref{temp2}. Since
$\lambda_{\min}\left((A^\star)^*A^\star+\gamma (A^*A-(A^\star)^*A^\star)\right)$ is a continuous function of $\gamma$, its supremum will be attained for some $\hat \gamma \in [0,1]$.
Thus by Lemma~\ref{l4} there exists an eigenvector $\hat v$ of $(A^\star)^*A^\star+\hat \gamma (A^*A-(A^\star)^*A^\star)$ corresponding to the optimal eigenvalue
$\lambda_{\min}\left((A^\star)^*A^\star+\hat \gamma (A^*A-(A^\star)^*A^\star)\right)$ such that
\begin{equation}\label{temp3}
\|\hat v\|=1\quad \text{and} \quad \hat v^*\left(A^*A-(A^\star)^*A^\star\right) \hat v=0.
\end{equation}
This implies that
\begin{eqnarray*}
2\,\lambda_{\min}\left((A^\star)^*A^\star+\hat \gamma (A^*A-(A^\star)^*A^\star)\right) &=& 2\,\hat v^*\left((A^\star)^*A^\star+\hat \gamma (A^*A-(A^\star)^*A^\star)\right)\hat v\\
&=& 2\,\hat v^* (A^\star)^*A^\star \hat v \quad (\because \text{from~\eqref{temp3}}) \\
&=& 2\,\max\left\{\frac{\|A\hat v\|^2}{\|\hat v\|^2},\frac{\|A^\star \hat v\|^2}{\|\hat v\|^2}\right\}\\
&\geq& \delta_0^{\text{pal}_\star}(A,A^\star)^2,
\end{eqnarray*}
where the last inequality follows from~\eqref{p1}. This shows equality in~\eqref{temp1} and hence the assertion.
\eproof

\subsection{Dissipative Hamiltonian structure}\label{subsec:dh}

In this section, we consider the  $n\times n$ DH pencils $s E+(J-R)$, where $E \succeq 0$, $R \succeq 0$, and $J^*=-J$. If $s E+(J-R)$ is  a real DH pencil, i.e., $J,R,Q \in \R^{n,n}$, then such a pencil is singular if and only if $E$ and $J-R$ have a common nonzero null space which is true if and only if there is a common nonzero vector in the null spaces of $J$, $R$, and $E$~\cite[Theorem 5]{MehMW20}. Thus the structured distance to singularity for a real DH pencil is equal to the structured distance to a common null space of $J$, $R$, and $E$.
This distance was considered and computable bounds were obtained in~\cite{MehMW20} for real DH pencils with respect to skew-symmetric perturbations to $J$, and symmetric but possibly indefinite perturbations to $R$ and $E$. We denote this distance from~\cite[Theorem 13]{MehMW20} by $\delta_0^i(J,R,E)$, where $i$ stands for the symmetric indefinite perturbations to $E$ and $R$. 

We study the distance to a common null space for complex DH pencils
$sE+(J-R)$ with respect to structure-preserving \emph{negative semidefinite} perturbations to $R$ and $E$, and skew-Hermitian perturbations to $J$.
In view of Lemma~\ref{lemdh1}, a DH pencil $s E+(J-R)$ has a common null space if and only if $J$, $R$, and $E$ have a common null space. Thus we consider the following distances.

\begin{definition}\label{def:nulldh}
Consider a regular complex DH pencils
$sE+(J-R)$, where $J,R,E \in \C^{n,n}$, $J^*=-J$, $R\succeq 0$, and 
$E\succeq 0$.
\begin{enumerate}
\item  The \emph{structured distance $\delta_0^d(J,R,E)$} to a common null space of $J$, $R$, and $E$ with respect to structure-preserving perturbations to J, and Hermitian negative semidefinite perturbations to $R$ and $E$ from the set 
\begin{equation}\label{eq:defst}
\mathcal S_d(R,E):= \left\{ (\Delta_E,\Delta_R)\in (\C^{n,n})^2:~\Delta_E\preceq 0,\, \Delta_R\preceq 0,\, E+\Delta_E \succeq 0,\,R+\Delta_R \succeq 0\right\}
\end{equation}
is defined by 
\begin{eqnarray}\label{eq:defjrq}
&\delta_0^d(J,R,E):=\inf\Big \{ \sqrt{{\|\Delta_J\|}^2+{\|\Delta_R\|}^2+{\|\Delta_E\|}^2}:~(\Delta_R,\Delta_E) \in \mathcal S_d(R,E),\, \Delta_J \in {\rm SHerm}(n)\nonumber \\
& \text{\rm ker}(E+\Delta_E) \cap \text{\rm ker}(R+\Delta_R) \cap \text{\rm ker}(J+\Delta_J) \neq \{0\}
\Big \}.
\end{eqnarray}

{\rm
Note that 
\begin{itemize}
\item the unstructured distance $\delta_0(J-R,E)$ when $J-R$ is treated as one matrix is obtained in Theorem~\ref{thm:unstrnull}. However, if  $J$ and $R$ are separately perturbed, then we consider the unstructured distance $\delta_0(J,R,E)$ which is defined when
 $\mathcal S_d(R,E)$ and ${\rm SHerm}(n)$ are respectively replaced by $(\C^{n,n})^2$ and $\C^{n,n}$ in~\eqref{eq:defjrq}.  Along the lines of Theorem~\ref{thm:unstrnull}, one can easily show that
 \begin{equation}\label{eq:unstr_DHnull}
 \delta_0(J,R,E)=\sqrt{\lambda_{\min}\left(-J^2+R^2+E^2\right)},
 \end{equation}
 and it is attained for rank-1 perturbations to $J$, $R$ and $E$;
\item 
$\delta_0^d(J,R,E) < \infty$ as the perturbations $(\Delta_J,\Delta_R,\Delta_E)=(-J,-R,-E)$ results in a zero pencil; 
\item  if $sE+(J-R)Q$ is a real DH pencil then by~\cite[Theorem 13]{MehMW20}, we have
\[
\delta_0^d(J,R,E) \geq \delta_0^i(J,R,E)  \geq  \delta_0(J,R,E).
\]
\end{itemize}

We also consider the structured distances to a common null space for DH pencils
when only any one/two matrices from $\{J,R,E\}$ are subject to perturbations.
}

\item The \emph{structured distance $\delta_0^d(J,R)$}, with respect to skew-Hermitian perturbations to $J$ and Hermitian negative semidefinite perturbations to $R$, is defined by
\begin{eqnarray}\label{def:nulljr}
&\delta_0^d(J,R):=\inf\Big \{ \sqrt{{\|\Delta_J\|}^2+{\|\Delta_R\|}^2}:~\Delta_R \in \C^{n,n},\,\Delta_R \preceq 0,\,(R+\Delta_R) \succeq 0,  \Delta_J \in {\rm SHerm}(n),\nonumber \\
& \text{\rm ker}(E) \cap \text{\rm ker}(R+\Delta_R) \cap \text{\rm ker}(J+\Delta_J) \neq \{0\}
\Big \}.
\end{eqnarray}
Similarly, we can define the \emph{structured radius $\delta_0^d(J,E)$}  with respect to  skew-Hermitian perturbations to $J$ and Hermitian negative semidefinite perturbations to $E$.

\item The \emph{structured distance $\delta_0^d(R,E)$}, with respect to Hermitian negative semidefinite perturbations to $R$ and $E$, is defined by
\begin{eqnarray*}
\delta_0^d(R,E):=\inf\Big \{ \sqrt{{\|\Delta_R\|}^2+{\|\Delta_E\|}^2}:~(R,E)\in \mathcal S_d(R,E),\nonumber  \text{\rm ker}(E+\Delta_E) \cap \text{\rm ker}(R+\Delta_R) \cap \text{\rm ker}(J) \neq \{0\}
\Big \}.
\end{eqnarray*}

\item The \emph{structured distance $\delta_0(J)$}, with respect to skew-Hermitian perturbations to $J$, is defined by
\begin{eqnarray*}
\delta_0(J):=\inf\Big \{ \|\Delta_J\|:~ \Delta_J \in {\rm SHerm}(n), \text{\rm ker}(E) \cap \text{\rm ker}(R) \cap \text{\rm ker}(J+\Delta_J) \neq \{0\}
\Big \}.
\end{eqnarray*}

\item The \emph{structured distance $\delta_0^d(R)$}, with respect to Hermitian negative semidefinite perturbations to $R$, is defined by
\begin{eqnarray*}
\delta_0^d(R):=\inf\Big \{ \|\Delta_R\|:~\Delta_R \in \C^{n,n},\,\Delta_R \preceq 0,\,(R+\Delta_R) \succeq 0, \text{\rm ker}(E) \cap \text{\rm ker}(R+\Delta_R) \cap \text{\rm ker}(J) \neq \{0\}
\Big \}.
\end{eqnarray*}
Similarly, we can define the \emph{structured distance $\delta_0^d(E)$}  with respect to   Hermitian negative semidefinite perturbations only to $E$.

\end{enumerate}
If the perturbations in $R$ and $E$ are restricted to be of rank one, then we denote the corresponding distance by adding an index $1$. For e.g., $\delta_{10}^d(J,R,E)$ denotes the structured distance to a common null space with respect to
 skew-Hermitian perturbations to $J$ and structure-preserving rank one perturbations to $R$ and $E$.
\end{definition}

In the following theorem, we obtain a formula for $\delta_0^d(J,R,E)$ in terms of sum of Rayleigh quotient and generalized Rayleigh quotients of some semidefinite matrices. 

\begin{theorem}\label{thm:mainJRE}
Let $sE+(J-R)$ be a regular DH pencil, where $J,R,E \in \C^{n,n}$ such that $J^*=-J$, $R\succeq 0$, and $E \succeq 0$. Then 
\[
\delta_0^d(J,R,E)^2=\min\left\{\delta_0^d(J,R)^2,\delta_0^d(J,E)^2,
\min_{\alpha \in \mathcal M_3} \frac{\alpha^*J^*J\alpha}{\alpha^*\alpha}+
 \left(\frac{\alpha^*R^2\alpha}{\alpha^*R\alpha}\right)^2+ \left(\frac{\alpha^*E^2\alpha}{\alpha^*E\alpha}\right)^2
\right\},
\]
where  $\mathcal M_3=\text{\rm ker}(R)^{c}\cap \text{\rm ker}(E)^{c}$, and $\delta_0^d(J,R)$ and $\delta_0^d(J,E)$ are given by Table~\ref{table:dhpencilformulas}.
\end{theorem}
\proof By definition~\eqref{eq:defjrq},
\begin{eqnarray}\label{eq:dhprf1}
\delta_0^d(J,R,E)~=&\inf\Big \{ \sqrt{{\|\Delta_J\|}^2+{\|\Delta_R\|}^2+{\|\Delta_E\|}^2}:~(\Delta_R,\Delta_E) \in \mathcal S_d(R,E),\, \Delta_J \in {\rm SHerm}(n), \nonumber \\
& \text{\rm ker}(E+\Delta_E) \cap \text{\rm ker}(R+\Delta_R) \cap \text{\rm ker}(J+\Delta_J) \neq \{0\}
\Big \}\nonumber \\
=&\inf\Big \{ \sqrt{{\|\Delta_J\|}^2+{\|\Delta_R\|}^2+{\|\Delta_E\|}^2}:~(\Delta_R,\Delta_E) \in \mathcal S_d(R,E),\, \Delta_J \in {\rm SHerm}(n),\nonumber \\
& x(\neq 0)\in \C^n,\,(J+\Delta_J)x=0,\,(R+\Delta_R)x=0,\,(E+\Delta_E)x=0
\Big \}\nonumber \\
=&\inf\Big \{ \sqrt{{\|\Delta_J\|}^2+{\|\Delta_R\|}^2+{\|\Delta_E\|}^2}:~(\Delta_R,\Delta_E) \in \mathcal S_d(R,E),\, \Delta_J \in {\rm SHerm}(n),\nonumber \\
& x(\neq 0)\in \C^n,\,\Delta_Jx=-Jx,\,\Delta_R x=-Rx,\,\Delta_Ex=-Ex
\Big \}.
\end{eqnarray}
Note that $\C^n=\mathcal M_1 \cup \mathcal M_2 \cup \mathcal M_3$, where
 $\mathcal M_1=\{x \in \C^n:~x^*Rx=0\}$,  $\mathcal M_2=\{x \in \C^n:~x^*Ex=0\}$, and  $\mathcal M_3=\{x \in \C^n:~x^*Rx\neq 0, \,x^*Ex \neq 0\}$. Since $R \succeq 0$, for any $x \in \C^n$ we have $x^*Rx \geq 0$, and $x^*Rx=0$  if and only if $Rx=0$. This implies that $\mathcal M_1=\text{\rm ker}(R)$, and similarly we have 
  $\mathcal M_2=\text{\rm ker}(E)$, and $\mathcal M_3=\text{\rm ker}(R)^c \cap \text{\rm ker}(E)^c$. This yields from~\eqref{eq:dhprf1} that 
 \begin{equation}\label{eq:minof3}
 \delta_0^d(J,R,E)^2=\min\{\mu_1,\mu_2,\mu_3\},
 \end{equation}
 where
\begin{eqnarray*}
\mu_1~=&\inf\Big \{ \sqrt{{\|\Delta_J\|}^2+{\|\Delta_E\|}^2}:~\Delta_E\in \C^{n,n},\,\Delta_E \preceq 0, E+\Delta_E \succeq 0,\, \Delta_J \in {\rm SHerm}(n),\nonumber \\
& x(\neq 0)\in \mathcal M_1,\,\Delta_Jx=-Jx,\,\Delta_Ex=-Ex
\Big \},
\end{eqnarray*}
\begin{eqnarray*}
\mu_2~=&\inf\Big \{ \sqrt{{\|\Delta_J\|}^2+{\|\Delta_R\|}^2}:~\Delta_R\in \C^{n,n},\,\Delta_R \preceq 0, R+\Delta_R \succeq 0,\, \Delta_J \in {\rm SHerm}(n),\nonumber \\
& x(\neq 0)\in \mathcal M_2,\,\Delta_Jx=-Jx,\,\Delta_Rx=-Rx
\Big \},
\end{eqnarray*}
and
\begin{eqnarray}\label{dfmu3}
\mu_3~=&\inf\Big \{ \sqrt{{\|\Delta_J\|}^2+{\|\Delta_R\|}^2+{\|\Delta_E\|}^2}:~(\Delta_R,\Delta_E) \in \mathcal S_d(R,E),\, \Delta_J \in {\rm SHerm}(n),\nonumber \\
& x(\neq 0)\in \mathcal M_3,\,\Delta_Jx=-Jx,\,\Delta_Rx=-Rx,\,\Delta_Ex=-Ex
\Big \}.
\end{eqnarray}
Observe from~\eqref{def:nulljr} that $\mu_1=\delta_0^d(J,E)$ and $\mu_2=\delta_0^d(J,R)$ which are obtained in Table~\ref{table:dhpencilformulas}. Thus it remains to compute $\mu_3$.
Since $ \mathcal S_d(R,E) \subseteq \{(\Delta_R,\Delta_E)\in (\C^{n,n})^2:~\Delta_R \preceq 0,\,\Delta_E \preceq 0\}$, we have from~\eqref{dfmu3} that
\begin{eqnarray}\label{dfmu3a}
&\mu_3 \geq \inf\Big \{ \sqrt{{\|\Delta_J\|}^2+{\|\Delta_R\|}^2+{\|\Delta_E\|}^2}:~\Delta_R,\Delta_E\in \C^{n,n},\,\Delta_R \preceq 0,\Delta_E \preceq 0 ,\, \Delta_J \in {\rm SHerm}(n),\nonumber \\
& x(\neq 0)\in \mathcal M_3,\,\Delta_Jx=-Jx,\,\Delta_Rx=-Rx,\,\Delta_Ex=-Ex
\Big \}.
\end{eqnarray}
If $\mathcal M_3 \neq \emptyset$, then the infimum in the right hand side of~\eqref{dfmu3a} is finite. Indeed, in view of Lemma~\ref{lemdh2}, for any $x(\neq 0)\in \mathcal M_3$ there exist $\Delta_R \preceq 0$ and $\Delta_E \preceq 0$ such that $\Delta_Rx=-Rx$ and $\Delta_Ex=-Ex$ if and only if $-x^*Rx <0$ and $-x^*Ex <0$. Clearly,  $-x^*Rx <0$ and $-x^*Ex <0$ since $x(\neq 0) \in \mathcal M_3$ implies that $Rx \neq 0$ and $Ex \neq 0$ which is true if and only if $x^*Rx >0$ and $x^*Ex >0$, since $E$ and $R$ are positive semidefinite. Also, from Lemma~\ref{l2} (since $\Delta \in {\rm Herm}(n)\Longleftrightarrow i\Delta \in {\rm SHerm}(n)$ ) there exists a skew-Hermitian $\Delta_J$ such that  $\Delta_Jx=-Jx$ if and only if $\real{x^*Jx}=0$. Clearly $\real{x^*Jx}=0$ since $J$ is skew-Hermitian. 

Next, we show that the inequality in~\eqref{dfmu3a} is actually an equality, since the infimum in the right hand side of~\eqref{dfmu3a}
is attained for some $\Delta_J,\Delta_R, \Delta_E$ such that $(\Delta_R,\Delta_E) \in  \mathcal S_d(R,E)$.
Indeed, for any $x \in \mathcal M_3$,  from Lemma~\ref{lemdh2} there exist 
\begin{equation}\label{eq:delrdelE}
\Delta_R=-\frac{(Rx)(Rx)^*}{x^*Rx} \quad \text{and}\quad 
\Delta_E=-\frac{(Ex)(Ex)^*}{x^*Ex}
 \end{equation}
 such that
 $ \Delta_R \preceq 0$, $\Delta_Rx=-Rx$, $ \Delta_E \preceq 0$, and $ \Delta_Ex=-Ex$. The mappings $\Delta_R$ and $\Delta_E$ are also of minimal spectral norm.
Also observe that $(R+\Delta_R)x=0$ and $(E+\Delta_E)x=0$. Thus in view of  Lemma~\ref{lemdh4}  we have $R+\Delta_R \succeq 0$ and $E+\Delta_E \succeq 0$. This implies that $(\Delta_R,\Delta_E) \in  \mathcal S_d(R,E)$.

Using $\Delta_R$ and $\Delta_E$ from~\eqref{eq:delrdelE},  and a minimal norm skew-Hermitian mapping from Lemma~\ref{l2},  we obtain 
\begin{eqnarray}\label{dfmu3b}
&\mu_3^2 = \inf\Big \{ {\|\Delta_J\|}^2+{\|\Delta_R\|}^2+{\|\Delta_E\|}^2:~\Delta_R,\Delta_E\in \C^{n,n},\,\Delta_R \preceq 0,\Delta_E \preceq 0 ,\, \Delta_J \in {\rm SHerm}(n),\nonumber \\
&x(\neq 0)\in \mathcal M_3,\,\Delta_Jx=-Jx,\,\Delta_Rx=-Rx,\,\Delta_Ex=-Ex
\Big \}\nonumber \\
& \hspace{-5cm}= \inf\Big \{\frac{{\|Jx\|}^2}{{\|x\|}^2}+\frac{{\|Rx\|}^4}{(x^*Rx)^2}+\frac{{\|Ex\|}^4}{(x^*Ex)^2}:~x(\neq 0)\in \mathcal M_3
\Big \},
\end{eqnarray}
which completes the proof. 
\eproof

We note that for a regular DH pencil $sE+(J-R)$, the other structured distances when the perturbations are restricted to only one/two matrices from $\{J,R,Q\}$,
can also be obtained by following the arguments similar to those of Theorem~\ref{thm:mainJRE}. We summarize these  results in Table~\ref{table:dhpencilformulas} and skip the details.

\begin{table}[!h]
\begin{center}
\caption{Various structured distances for a regular DH pencil $sE+(J-R)$.}
\label{table:dhpencilformulas}
\begin{tabular}{ |c|c|l| }
\hline
\multicolumn{3}{ |c| }{Columns of $U$ form an orthonormal basis for $\Omega$.} 
\\
\hline
\multicolumn{3}{ |c| }{$L_1$ and $L_2$ are respectively the Cholesky factors of $U^*RU$ and $U^*EU$.} 
\\
\hline
\multicolumn{3}{ |c| }{$f_\alpha(J,R)=\frac{\alpha^*J^*J\alpha}{\alpha^*\alpha}+(\frac{\alpha^*R^2\alpha}{\alpha^*R\alpha})^2$, $g_\alpha(R,E)=(\frac{\alpha^*R^2\alpha}{\alpha^*R\alpha})^2+(\frac{\alpha^*E^2\alpha}{\alpha^*E\alpha})^2$.} 
\\
\hline 
\hline
  \textbf{ Distance} & \textbf{$\Omega$} & \multicolumn{1}{ |c| }{\textbf{ Value}} \\ \hline
  \multirow{2}{*}{$\delta_0(J)^2$} & \multirow{2}{*}{$\text{\rm ker}(R)\cap \text{\rm ker}(E)$} & $\infty$ \quad \quad if $\Omega=\{0\}$ \\\cline{3-3}
& & $\left(\sigma_{\min}(JU)\right)^2$,\quad otherwise\\ \hline
\multirow{2}{*}{$\delta_0^d(R)^2$} & \multirow{2}{*}{$\text{\rm ker}(J)\cap \text{\rm ker}(E)$} & $\infty$ \quad \quad if $\Omega=\{0\}$ \\\cline{3-3}
& & $\left(\lambda_{\min}(L_1^{-*}U^*R^2UL_1^{-1})\right)^2$,\quad otherwise\\ \hline
\multirow{2}{*}{$\delta_0^d(E)^2$} & \multirow{2}{*}{$\text{\rm ker}(J)\cap \text{\rm ker}(R)$} & $\infty$ \quad \quad if $\Omega=\{0\}$ \\\cline{3-3}
& & $\left(\lambda_{\min}(L_2^{-*}U^*E^2UL_2^{-1})\right)^2$,\quad otherwise\\ \hline
\multirow{2}{*}{$\delta_0^d(J,R)^2$} & \multirow{2}{*}{$\text{\rm ker}(E)\cap \text{\rm ker}(R)^{c} $} & $\infty$ \quad \quad if $\text{\rm ker}(E)=\{0\}$ \\\cline{3-3}
& & $\min\{\delta_0(J)^2,\inf_{\alpha \in \Omega }f_\alpha(J,R)\}$,\quad otherwise\\ \hline
\multirow{2}{*}{$\delta_0^d(J,E)^2$} & \multirow{2}{*}{$\text{\rm ker}(R)\cap \text{\rm ker}(E)^{c} $} & $\infty$ \quad \quad if $\text{\rm ker}(R)=\{0\}$ \\\cline{3-3}
& & $\min\{\delta_0(J)^2,\inf_{\alpha \in \Omega}f_\alpha(J,E)\}$,\quad otherwise\\ \hline
\multirow{2}{*}{$\delta_0^d(R,E)^2$} & \multirow{2}{*}{$\text{\rm ker}(J)\cap\text{\rm ker}(R)^c\cap \text{\rm ker}(E)^{c} $} & $\infty$ \quad \quad if $\text{\rm ker}(J)=\{0\}$ \\\cline{3-3}
& & $\min\{\delta_0^d(R)^2,\delta_0^d(E)^2,\inf_{\alpha \in \Omega}g_\alpha(R,E)\}$,\quad otherwise\\ \hline
\end{tabular}
\end{center}
\end{table}

\begin{remark}\label{rem:rankonedh}{\rm 
Observe that (i)  any Hermitian matrix of rank one is necessarily semidefinite and (ii) only negative semidefinite perturbations
$\Delta_R$ and $\Delta_E$ in $R+\Delta_R \succeq 0$ and $E+\Delta_E \succeq 0$ can make a nonzero vector in the common null space of $R+\Delta_R$ and $E+\Delta_E$. In view of this, we see that for a skew-Hermitian perturbation $\Delta_J$ to $J$ and any structure-preserving rank one perturbations $\Delta_R$ and $\Delta_E$ to $R$
and $E$ such that ${\rm ker}(J) \cap {\rm ker}(R) \cap {\rm ker}(E) \neq \{0\}$ imply that $\Delta_R$ and $\Delta_E$ are necessarily negative semidefinite and 
$\sqrt{{\|\Delta_J\|}^2+{\|\Delta_R\|}^2+{\|\Delta_E\|}^2} \geq \delta_0^d(J,R,E)$. Also the minimal norm mappings $\Delta_R$ and $\Delta_E$ in the proof of Theorem~\ref{thm:mainJRE} can be chosen to be of rank one. Consequently, we have 
\[
\delta_{10}^d(J,R,E)=\delta_0^d(J,R,E)=\delta_{10}^i(J,R,E).
\] 

%

}
\end{remark}

\section{Structured distance to a common null space for higher degree polynomials}\label{sec:highdeg}

In this section, we show that the linear algebra approach proposed in 
Section~\ref{sec:dist_null}
can be generalized to matrix polynomials of the form $P(s)=\sum_{j=0}^m s^jA_j$, where
$(A_0,\ldots,A_m) \in \mathbb S \subseteq (\C^{n,n})^{m+1}$.
%
The structured distance to a common null space for $P(s)$ is denoted by
$\delta_0^{\mathbb S}(P(s))$ and defined as
\begin{eqnarray}
\delta_0^{\mathbb S}(P(s)):= \inf\left \{
\sqrt{\sum_{j=0}^m {\|\Delta_j\|}^2}:~(\Delta_0,\ldots,\Delta_m) \in\mathbb S,~\bigcap_{j=0}^m \text{\rm ker}(A_j-\Delta_j) \neq \{0\}
\right \}.
\end{eqnarray}
For the unstructured distance $\delta_0(P(s)):=\delta_0^{\mathbb S}(P(s))$ (when $\mathbb S = (\C^{n,n})^{m+1}$), and for the polynomials with  Hermitian and related structures, the results of Sections~\ref{subsec:unstrnullspace} and~\ref{subsec:null_herm} are easily extendable. We summarized the results for these cases in Table~\ref{tab:dist_null_poly} and skip the details.

\begin{table}[!h]
\caption{Structured distance to a common null space for matrix polynomials} \label{tab:dist_null_poly}
\begin{center}
\begin{tabular}{|c|c|c|}
\hline
structure & property & distance to a common null space \\ \hline
unstructure poly. & $\mathbb S = (\C^{n,n})^{m+1}$ & $\delta_0(P(s))=\sqrt{\lambda_{\min}(\sum_{j=0}^mA_j^* A_j)}$ \\ \hline
Hermitian & $A_j^*=A_j$ & $\delta_0^{{\rm Herm}}(P(s))=\delta_0(P(s))$ \\ \hline
skew-Hermitian & $A_j^*=-A_j$ & $\delta_0^{{\rm SHerm}}(P(s))=\delta_0^{{\rm Herm}}(iP(s))$ \\ \hline
$*$-even & $A_j^*=(-1)^jA_j^*$ & $\delta_0^{\mathbb S_e}(P(s))=\delta_0^{{\rm Herm}}(P(is))$ \\ \hline
$*$-odd & $A_j^*=(-1)^{j+1}A_j^*$ & $\delta_0^{\mathbb S_o}(P(s))=\delta_0^{{\rm Herm}}(iP(is))$ \\ \hline
\end{tabular}
\end{center}
\end{table}

Although, the ideas of Section~\ref{sec:pal_null} for the distance to a common null space for $\star$-palindromic polynomial
$P(s)=\sum_{j=0}^m s^jA_j$, where
$A_j={A}_{m-j}^\star$, and $\star=*$ for
$*$-palindromic structure and $\star=T$ for $T$-palindromic structure, can also be generalized to some extent but the polynomial case needs to be stated separately. The proof of Theorem~\ref{thm:res_null_pal} uses Lemma~\ref{l4}, which is no longer valid if the eigenvalue function is depending on more than one parameters. Instead, we use the following lemma which is a weaker version of Lemma~\ref{l4} to tackle the polynomial case. 

\begin{lemma}{\rm \cite{Bau85}}\label{lem:weakdifcond}
Let $G,H \in {\rm Herm}(n)$ and let the map $L:\R \rightarrow \R$ be given by $L(t):=\lambda_{\min}(G+tH)$. If $\lambda_{\min}(G)$ is a simple eigenvalue of $G$, then $L$ is differentiable at $t=0$. Moreover, there exists $u \in \mathcal{E}_{\lambda_{\min}}(G)$ such that  $\|u\|=1$ and 
$\frac{d}{dt}L(0)=u^*Hu$.
\end{lemma}

The Lemma~\ref{lem:weakdifcond} allows us to obtain the distance to a common null space for $\star$-palindromic polynomials of higher degree. 
\begin{theorem}
Let $P(s)=\sum_{j=0}^m s^jA_j$ be $\star$-palindromic and let $k=\lfloor{\frac{m-1}{2}}\rfloor$. Then
\begin{equation*}
\hat \lambda_{\min}=\sup_{\gamma_0,\ldots,\gamma_k \in [0,1]} \lambda_{\min}(f(\gamma_0,\ldots,\gamma_k)),
\end{equation*}
where $f(\gamma_0,\ldots,\gamma_k)=\sum_{j=0}^{k} (A_j^\star)^*A_j^\star+\gamma_j (A^*A-(A^\star)^*A^\star)$ if $m$ is odd, and
$f(\gamma_0,\ldots,\gamma_k)=\sum_{j=0}^{k} (A_j^\star)^*A_j^\star+\gamma_j (A^*A-(A^\star)^*A^\star+A^*_{\frac{m}{2}}A_{\frac{m}{2}})$ if $m$ is even, is attained for some $\hat \gamma_0,\ldots,\hat \gamma_k \in [0,1]$. Furthermore,
\begin{equation}\label{eq:ineq_null_pal}
\delta_0^{{\rm Pal}_\star}(P(s)) \geq 2 \,\hat \lambda_{\min},
\end{equation}
and equality holds in~\eqref{eq:ineq_null_pal} if $m \leq 2$ or $\hat \lambda_{\min}$ is a simple eigenvalue of $f(\hat \gamma_0,\ldots,\hat \gamma_k)$.
\end{theorem}
 \proof
The proof follows by (i) generalizing the steps of Theorem~\ref{thm:res_null_pal} for the $\star$-palindromic polynomials and (ii) in view of Lemma~\ref{lem:weakdifcond} using the fact that $\lambda_{\min}(f(\gamma_0,\ldots,\gamma_k))$
is differentiable at $\hat \gamma_0,\ldots,\hat \gamma_k$ when
$\hat \lambda_{\min}$ is a simple eigenvalue of $f(\gamma_0,\ldots,\gamma_k)$.
 \eproof
 
 We note that the results for DH pencils obtained in Section~\ref{subsec:dh} can be generalized in a straightforward way for matrix polynomials with DH related structure of the form $P(s)=-s^j J+ \sum_{i=0}^n s^iA_i$, where $j,n \geq 0$, $J,A_i \in \C^{n,n}$ such that $J^*=-J$ and $A_i \succeq 0$ for each $i=0,\ldots,n$. Such polynomials were considered in~\cite{MehMW20}. 

\section{Structured distance to singularity}\label{sec:dist_sing}

Note that  $\delta^{\mathbb S}(A,E) \leq \delta_0^{\mathbb S}(A,E)$, as $\delta_0^{\mathbb S}(A,E)$ is the distance to singularity through the common null space of the matrix pencil $L(s)=A+sE$.
This gives an upper bound
for the distance to singularity as in general such a situation need not occur for singular pencils.

In this section, we obtain a family of lower bounds for the structured distance to singularity $\delta^\mathbb S(A,E)$ in terms of structured eigenvalue backward errors. The structured eigenvalue backward error of  a $\lambda \in \C$ as an approximate eigenvalue of the pencil $L(s)=A+sE$ is denoted by $\eta^\mathbb S(A,E,\lambda)$ and defined as
\begin{eqnarray}\label{def:eigerror}
\eta^\mathbb S(A,E,\lambda):=\inf\big\{ \sqrt{\|\Delta_A\|^2+\|\Delta_E\|^2}\;:\;(\Delta_A,\Delta_E)\in \mathbb S,\,
\text{\rm det}((A-\Delta_A)+\lambda(E-\Delta_E))=0
\bigg\}.
\end{eqnarray}

When $S=(\C^{n,n})^2$, $\eta(A,E,\lambda):=\eta^\mathbb S(A,E,\lambda)$ is called the unstructured eigenvalue backward error. 

\paragraph{\textbf{First lower bound}}:~Let $(\Delta_A,\Delta_E) \in \mathbb S$ be such that $(A-\Delta_A)+s(E-\Delta_E)$ is singular, then $\text{\rm det}((A-\Delta_A)+\lambda(E-\Delta_E))=0$ for all $\lambda \in \C$. This implies in view of~\eqref{def:eigerror} that for any $\lambda_0 \in \C$, we have 
$\sqrt{\|\Delta_A\|^2+\|\Delta_E\|^2} \geq \eta^\mathbb S(A,E,\lambda)$. This implies that  $\sqrt{\|\Delta_A\|^2+\|\Delta_E\|^2} \geq \sup_{\lambda_0 \in \C}\eta^\mathbb S(A,E,\lambda_0)$. This yields a lower bound for $\delta^\mathbb S(A,E)$ as 
\begin{equation}\label{eq:lboundgenlam}
\delta^\mathbb S(A,E) \geq \sup_{\lambda \in \C}\eta^\mathbb S(A,E,\lambda).
\end{equation}

\paragraph{\textbf{A family of lower bounds}}:~Also note that $A+sE$ is singular if and only if there exist distinct $\lambda_1,\ldots,\lambda_d \in \C$ , $d\geq n+1$ such that
$\text{\rm det}(A+\lambda_i E)=0$ for all $i=1,\ldots,d$. This  
reformulates the distance to singularity problem into an equivalent problem of computing the nearest pencil with any pre-specified  $n+1$ distinct complex numbers $\lambda_1,\ldots,\lambda_{n+1}$ as its eigenvalues. Further note that this reformulation is independent of the choice of $n+1$ distinct complex numbers $\lambda_1,\ldots,\lambda_{n+1}$. This results in a family of lower bounds:
\begin{equation}\label{eq:famlbounds}
\delta^\mathbb S(A,E) \geq \max_{i=1,\ldots,n+1}\eta^\mathbb S(A,E,\lambda_i).
\end{equation}

Explicit formulae for the eigenvalue backward error of matrix polynomials have been derived in~\cite{BorKMS14,BorKMS15} for the structures listed in Table~\ref{tab:dist_null_poly} and palindromic polynomials. We recall these results from~\cite{BorKMS14,BorKMS15} for the pencil case and state them here in the form that allows us to write explicit lower bounds for $\delta^\mathbb S(A,E)$. For this, let us define
\begin{equation}\label{Eq:indchoice}
\mathcal K:= \left\{
(\lambda_1,\ldots,\lambda_{n+1})\in \C^{n+1}\;:\; \lambda_i \neq \lambda_j\,\text{ for}\,
 i \neq j,\; L(\lambda_i) \,\text{is invertible for all}\,i=1,\ldots,n+1
\right\}.
\end{equation}

\subsection{Unstructured distance to singularity}\label{subsec:unstrdistsing}

Suppose $\mathbb S=(\C^{n,n})^2$, then from~\eqref{eq:famlbounds} we have that 
\begin{equation}\label{eq:famlboundsuns}
\delta(A,E) \geq \max_{1=1,\ldots,n+1}\eta(A,E,\lambda_i),
\end{equation}
where $\delta(A,E)$ and $\eta(A,E,\lambda_i)$ are the unstructured distance to singularity and the unstructured eigenvalue backward error defined by~\eqref{sdistance} and~\eqref{def:eigerror}, respectively. 
A direct application of~\cite[Theorem 4.1]{BorKMS14} in~\eqref{eq:famlboundsuns} gives a lower bound for the unstructured distance to singularity $\delta(A,E)$, as stated in the following theorem. 
\begin{theorem}\label{thm:bacerro}
Let $L(s)=A+sE$ be a regular pencil, where $A,E \in \C^{n,n}$, and let $(\lambda_1,\ldots,\lambda_{n+1}) \in \mathcal K$. Then
\begin{equation}
\delta(A,E)^2 \geq
\max_{i=1,\ldots,n+1} \frac{1}{\lambda_{\max}(H_i)},
\end{equation}
where $H_i=\mat{c}1 \\ \bar \lambda_i \rix M_i^*M_i \mat{cc}1 & \lambda_i\rix$ and $M_{i}=(L(\lambda_i))^{-1}$ for each $i=1,\ldots,n+1$.
\end{theorem}

Next, we state a result
similar to Theorem~\ref{thm:bacerro} without its proof to give lower bounds for $\delta(A)$ and $\delta(E)$.
\begin{theorem}
Let $L(s)$ be a regular matrix pencil, where $A,E \in \C^{n,n}$, and let
$(\lambda_1,\ldots,\lambda_{n+1})\in \mathcal K$. Then
\begin{enumerate}
\item $\delta(A)^2 \geq  \max_{i=1,\ldots,n+1}\frac{1}{\lambda_{\max}(H_i)}$, and
\item $\delta(E)^2 \geq   \max_{i=1,\ldots,n+1}\frac{1}{|\lambda_i|^2\lambda_{\max}(H_i)}$,
\end{enumerate}
where $H_i=M_i^*M_i$ and $M_i=(L(\lambda_i))^{-1}$ for each $i=1,\ldots,n+1$.
\end{theorem}

\subsection{Hermitian and related structure}

Here, we consider Hermitian pencils $L(s)=A+sE$, where $(A,E) \in \mathbb S=(\text{Herm}(n))^2$. Then from~\eqref{eq:famlbounds} we have that 
\begin{equation}\label{eq:famlboundsherm}
\delta^{{\rm Herm}}(A,E) \geq \max_{i=1,\ldots,n+1}\eta^{\rm{Herm}}(A,E,\lambda_i),
\end{equation}
where $\delta^{{\rm Herm}}(A,E)$ and $\eta^{\rm{Herm}}(A,E,\lambda_i)$ are the Hermitian distance to singularity and the Hermitian eigenvalue backward error defined by~\eqref{sdistance} and~\eqref{def:eigerror}, respectively. 
As an application of~\cite[Theorem 4.3]{BorKMS14} in~\eqref{eq:famlboundsherm} yields a lower bound to the structured distance to singularity for Hermitian pencils. More precisely, we have the following result.
\begin{theorem}\label{thm:lbounddistrherm}
Let $L(s)=A+sE$ be a regular Hermitian pencil, where $A,E \in \text{Herm}(n)$, and let $(\lambda_1,\ldots,\lambda_{n+1}) \in \mathcal K\setminus \R^{n+1}$. Then
\[
\delta^{\rm{Herm}}(A,E)^2\geq  \max_{j=1,\ldots,n+1}
 \Big(\min_{t_0,t_1\in \R}\lambda_{\max}(G_j+t_0H_{1j}+t_1H_{2j})\Big)^{-1},
\]
where for each $j=1,\ldots,n+1$
\[
G_j=\mat{c}1 \\ \bar \lambda_j \rix M_j^*M_j \mat{cc}1 & \lambda_j\rix,
\quad H_{1j}=i\begin{bmatrix}  M_j - M_j^* & \lambda_j M_j\\ -\bar{\lambda}_jM_j^* & 0 \end{bmatrix},\quad H_{2j} = i\begin{bmatrix} 0 & -M_j^*\\ M_j & \lambda_j M_j-\bar{\lambda}_jM_j^* \end{bmatrix},
\]
and $M_j=(L(\lambda_j))^{-1}$.
\end{theorem}

\begin{remark}\label{rem:hermrelated}{\rm
In view of Table~\ref{tab:dist_sing_herrel},  a family of lower bounds for the structured distance to singularity for pencils
with skew-Hermtian, $*$-even, and $*$-odd structures can also be obtained using Theorem~\ref{thm:lbounddistrherm}.  In Table~\ref{tab:dist_sing_herrel}, $\eta^{\text{SHerm}}(A,E,\lambda)$, $\eta^{\rm{\mathbb S_e}}(A,E,\lambda)$, and $\eta^{\rm{\mathbb S_o}}(A,E,\lambda)$ respectively denote the skew-Hermitian, $*$-even, and $*$-odd eigenvalue backward errors.

\begin{table}[!h]
\begin{center}
\caption{Structured eigenvalue backward errors for Hermitian related structures} \label{tab:dist_sing_herrel}
\begin{tabular}{|c|c|}
\hline
structure & relation with Hermitian distance \\ \hline
skew-Hermitian & $\eta^{\text{SHerm}}(A,E,\lambda)=\eta^{\text{Herm}}(iA,iE,\lambda)$\\ \hline
$*$-even & $\eta^{\rm{\mathbb S_e}}(A,E,\lambda)=\eta^{\text{Herm}}(A,iE,-i\lambda) $\\ \hline
$*$-odd & $\eta^{\rm{\mathbb S_o}}(A,E,\lambda)=\eta^{\text{Herm}}(iA,iE,-i\lambda)$ \\ \hline
\end{tabular}
\end{center}
\end{table}
}
\end{remark}

Along the lines of Theorem~\ref{thm:lbounddistrherm}, we obtain lower bounds for $\delta^{{\rm Herm}}(A)$ and $\delta^{{\rm Herm}}(E)$. 
An analogous result for $\delta^{{\rm Herm}}(A)$ and $\delta^{{\rm Herm}}(E)$ without its proof is stated in the following.

\begin{theorem}
Let $L(s)=A+sE$ be a regular Hermitian pencil, where $A,E \in {\rm Herm}(n)$ and let $(\lambda_1,\ldots,\lambda_{n+1}) \in \mathcal K\setminus \R^{n+1}$. Then
\begin{enumerate}
\item $\delta^{{\rm Herm}}(E)=\infty$ if $A$ is invertible, and
\[
\delta^{{\rm Herm}}(E) \geq \max_{j=1,\ldots,n+1}
\left(\min_{t\in \R} \lambda_{\max}\left(\widetilde G_j+t \widetilde H_j\right)
\right)^{-1},
\]
otherwise, where for each $j=1,\ldots,n+1$,  $\widetilde G_j=|\lambda_j|^2M_j^*M_j$, $\widetilde H_j=i(\lambda_j M_j-\bar \lambda_j M_j^*)$, and $M_j=(L(\lambda_j))^{-1}$;
\item we have
\[
\delta^{{\rm Herm}}(A) \geq  \max_{j=1,\ldots,n+1}
\left(\min_{t\in \R} \lambda_{\max}\left(G_j+tH_j\right)
\right)^{-1},
\]
where for each $j=1,\ldots,n+1$,  $G_j=M_j^*M_j$, $H_j=i(M_j-M_j^*)$, and $M_j=(L(\lambda_j))^{-1}$.
\end{enumerate}
\end{theorem}
%

\subsection{Palindromic structure}
In this section, we obtain lower bounds for $\star$-palindromic pencil $L(s)=A+sE$, where $E=A^\star$, and $\star=*$ if $L(s)$ is $*$-palindromic and $\star=T$ if $L(s)$ is T-palindromic. 
From~\eqref{eq:famlbounds} we have that 
\begin{equation}\label{eq:famlboundspal}
\delta^{{\rm Pal}_\star}(A,A^\star) \geq \max_{i=1,\ldots,n+1}\eta^{\rm{Pal}_\star}(A,A^\star,\lambda_i),
\end{equation}
where $\delta^{{\rm Pal}_\star}(A,A^\star)$ and $\eta^{\rm{Pal}_\star}(A,A^\star,\lambda_i)$ are the $\star$-palindromic distance to singularity and the $\star$-palindromic eigenvalue backward error defined by~\eqref{sdistance} and~\eqref{def:eigerror}, respectively. 
By applying~\cite[Theorem 4.6]{BorKMS15} for $*$-palindromic pencils and~\cite[Theorem 5.2]{BorKMS15} for T-palindromic pencils in~\eqref{eq:famlboundspal}, we obtain a  lower bound for  the structured distance to singularity $\delta^{{\rm pal}_{\star}}(A,A^\star)$. More precisely, we have the following result.

\begin{theorem}\label{thm:lboundpalpen}
Let $L(s)=A+sA^\star$ be a regular $\star$-palindromic pencil, where $A\in \C^{n,n}$ and let $(\lambda_1,\ldots,\lambda_{n+1}) \in \mathcal K$. For each $\lambda_j$, define $M_j=L(\lambda_j)^{-1}$, 
$C_j:=\mat{cc} M_j^\star & 0 \\ \bar \lambda_j M_j^\star-M_j & -\lambda_j M_j \rix$, $\Gamma_j:=\mat{cc}\sqrt{\frac{1}{1+|\lambda_j|}}I_n &0\\
0& \sqrt{\frac{|\lambda_j|}{1+|\lambda_j|}}I_n\rix$, and
$G_j:=\Gamma_j^{-1}\mat{c}I_n \\ \bar \lambda_j I_n \rix M_j^*M_j \mat{cc}I_n&  \lambda_j I_n \rix \Gamma_j^{-1}$. Then
\begin{enumerate}
\item when $\star=*$, we have
\begin{equation}\label{resultpal}
\delta^{\rm{pal_*}}(A,A^*)^2 \geq 2
\max_{j=1,\ldots,n+1}
\left(\min_{t_1,t_2\in \R}
\lambda_{\max}(G_j+t_1H_{1j}+t_2H_{2j})
\right)^{-1},
\end{equation}
where $H_{1j}:=\Gamma_j^{-1}(C_j+C_j^*)\Gamma_j^{-1}$ and  $H_{2j}:=i\Gamma_j^{-1}(C_j-C_j^*)\Gamma_j^{-1}$.
\item when $\star=T$, we have
\begin{equation}\label{resultpalT}
\delta^{\rm{pal_T}}(A,E)^2 \geq 2
\max_{j=1,\ldots,n+1}
\left(\min_{t\in [0,\infty)}
\lambda_{2}\left(\mat{cc}G_j &t\overline {S_j} \\ tS_j & G_j
\rix\right)
\right)^{-1},
\end{equation}
where $S_j=\Gamma_j^{-1}\left(C_j+C_j^T\right)\Gamma_j^{-1}$.
\end{enumerate}
\end{theorem}

\section{Unstructured vs. structured distances}\label{sec:numresults}

In this section, we illustrate the significance of our structured distances on some randomly generated structured pencils and compare them with the unstructured ones.
As far as we know, no other work has been done in the literature on  the structured distance to singularity except the one in~\cite{MehMW20} for DH pencils. Therefore, we cannot compare to existing work. 

To compute the distances in all cases, we used the \emph{GlobalSearch} in Matlab Version No. 9.5.0 (R2018b) to solve the associated optimization problem, except for the computation of lower bounds to $\delta^{{\rm Herm}}(A,E)$ and  $\delta^{{\rm pal}_*}(A,E)$ where we first  used the software package CVX~\cite{GraBY08}
to solve the inner optimization and then the  \emph{GlobalSearch} for outer optimization. In the following \emph{l.b.}  stands for the term ``lower bound" and 
\emph{u.b.} stands for the term ``upper bound".

\begin{example}{\rm (For Hermitian pencils)
~In the first experiment, we generate random Hermitian pencils $A+sE$, where $A,E \in {\rm Herm}(n)$ for different values of $n$ and record the results in Table~\ref{Tab:Hermresult}. The second column records the lower bound to unstructured distance  $\delta(A,E)$ obtained in Theorem~\ref{thm:bacerro} and the third column records the lower bound to  structured distance $\delta^{{\rm Herm}}(A,E)$ obtained in Theorem~\ref{thm:lbounddistrherm}. As expected, two lower bounds are  significantly different. 
This shows that the Hermitian pencil is more robustly regular under the 
structure-preserving perturbations. Note that for Hermitian pencils the unstructured and the structured distance to a common null space are equal (see Theorem~\ref{thm:hermnull}) which is an upper bound to $\delta^{{\rm Herm}}(A,E)$. This shows that the lower bound obtained in this paper gives a good estimation of the actual structured distance to singularity. 

\begin{center}
\begin{table}[h!]
\begin{center}
\caption{Various structured distances for Hermitian  pencils}
\label{Tab:Hermresult}
\begin{tabular}{|c|c|c|c|}
\hline
size  & \emph{l.b.} to $\delta(A,E)$  &\emph{l.b.} to $\delta^{{\rm Herm}}(A,E)$ &$\delta_0(A,E) =\delta_0^{{\rm Herm}}(A,E)$\\ 
& (Theorem~\ref{thm:bacerro}) &(Theorem~\ref{thm:lbounddistrherm}) & (Theorem~\ref{thm:hermnull}) \\ \hline
$3 \times 3$ &1.3922 &   1.5725  &  2.0344 \\ \hline
$4 \times 4$ & 2.2381  &   2.9006  &  3.8731 \\ \hline
$5 \times 5$ &2.1051  &  2.8435&    2.8968 \\ \hline
$6 \times 6$ &1.3709   & 2.1892    &2.3438 \\ \hline
$7 \times 7$ & 1.3560   & 2.5681   & 3.0594 \\ \hline
$8 \times 8$  & 1.6515   & 2.6032   & 3.0628 \\ \hline
\end{tabular}
\end{center}
\end{table}
\end{center}

}
\end{example}

\begin{example}{\rm  {(For Palindromic pencils)}
~In the second experiment, we generate random $*$-palindromic pencils $A+sA^*$, where $A \in \C^{n,n}$ for different values of $n$, and record the results in Table~\ref{Tab:Palresults}. We get similar results as for Hermitian pencils. 
The $*$-palindromic pencils are more robustly regular under structure preserving perturbations. As proven in Theorem~\ref{thm:res_null_pal}, the structured distance $\delta_0^{{\rm pal}_*}(A,A^*)$ is different than the unstructured distance $\delta_0(A,A^*)$, this fact is also illustrated numerically in Table~\ref{Tab:Palresults}.

\begin{center}
\begin{table}[h!]
\begin{center}
\caption{Various structured distances for $*$-palindromic  pencils}
\label{Tab:Palresults}
\begin{tabular}{|c|c|c|c|c|}
\hline
size  & \emph{l.b.} to $\delta(A,A^*)$  &\emph{l.b.} to $\delta^{{\rm pal}_*}(A,A^*)$ &$\delta_0(A,A^*)$ 
& $\delta_0^{{\rm pal}_*}(A,A^*)$\\ 
& (Theorem~\ref{thm:bacerro}) &(Theorem~\ref{thm:lboundpalpen}) & (Theorem~\ref{thm:unstrnull}) &(Theorem~\ref{thm:res_null_pal}) \\ \hline
$3 \times 3$&  1.2011   & 1.5953  &  1.9699 &   2.1048 \\ \hline
$4\times 4$ & 1.1254    &1.4374 &   2.2993 &   2.3528 \\ \hline
$5 \times 5$ & 0.9150  &  1.0071 &   2.0903 &   2.1146 \\ \hline
$6 \times 6$ & 1.0142   & 1.2130 &   1.7797 &   1.7815 \\ \hline
$7 \times 7$ & 0.8964    &1.0243 &   2.1385 &   2.1578 \\ \hline
$8 \times 8$ & 0.2175   & 0.6294 &   1.5865 &   1.6430 \\ \hline
\end{tabular}
\end{center}
\end{table}
\end{center}

}
\end{example}

\begin{example}{\rm (DH pencils)
~In our third experiment, we consider DH pencils and compare our results with the ones in~\cite{MehMW20}. In our notation $\delta_0^i(J,R,E)$ which is also equal to the structured distance to singularity for DH pencils, is obtained in~\cite[Theorem 13]{MehMW20}, where the perturbations in $J$, $R$, and $E$ were measured with respect to the Frobenius norm ${\|[\Delta_J,\Delta_R,\Delta_E]\|}_F$. Clearly, the unstructured distance $\delta_0(J,R,E)=\sqrt{\lambda_{\min}(-J^2+R^2+E^2)}$ from~\eqref{eq:unstr_DHnull}
gives a lower bound to  $\delta_0^i(J,R,E)$ which coincides with the lower bound in~\cite[Theorem 13]{MehMW20}. An upper bound to $\delta_0^i(J,R,E)$ was given by
$\sqrt{2}\delta_0(J,R,E)$ in~\cite[Theorem 13]{MehMW20}.

The structured distance $\delta_0^d(J,R,E)$ with respect to skew-Hermitian perturbations to $J$ and negative semidefinite perturbations to $R$ and $E$ also
gives an upper bound to $\delta_0^i(J,R,E)$. To compare it with the upper bound
$\sqrt{2}\delta_0(J,R,E)$ in~\cite{MehMW20}, we compute the Frobenius norm
${\|[\Delta_J,\Delta_R,\Delta_E]\|}_F$ of our optimal perturbation for which 
$\delta_0^d(J,R,E)$ was attained.

In Table~\ref{Tab1:DHresults}, we record $\delta_0(J,R,E)$, $\delta_0^d(J,R,E)$, and two upper bounds (the better of two is in bold) for DH pencils $sE-(J-R)$, where
$J=\mat{cc}0 & -0.5\\ 0.5 & 0\rix$
and $R=\mat{cc}0.18 &  0.42\\0.42 & 1.03\rix$ are fixed, and $E$ is chosen differently. We note that (i) when $E=\text{diag}(0,1)$~\cite[Example 24]{MehMW20}, we have $0.5819 \leq \delta_0^i(J,R,E) \leq 0.8229$ which is 
significantly smaller than the structured distance $\delta_0^d(J,R,E)=1.1181$.
This shows that structured distance with respect to semidefinite perturbations to $R$ and $E$ is more robust. In this case, the upper bound $0.8229$ by~\cite{MehMW20} is better than ours; (ii) for the other two pencils when
$E=\text{diag}(1,0)$ or $E=\text{diag}(1,1)$, our upper bound improves the upper bound in~\cite{MehMW20}.

\begin{center}
\begin{table}[h!]
\begin{center}
\caption{structured distances for DH pencils inspired by~\cite[Example 24]{MehMW20}}
\label{Tab1:DHresults}
\begin{tabular}{|c|c|c|c|c|}
\hline 
$E$ &   $\delta_0(J,R,E)$&  $\emph{u.b.} (=\sqrt{2}\delta_0(J,R,E)$)& our \emph{u.b.} &  $\delta_0^d(J,R,E)$\\ 
&\eqref{eq:unstr_DHnull} & \cite[Theorem 13]{MehMW20} &  & Theorem~\ref{thm:mainJRE} \\
\hline
 $\text{diag}(0,1)$&0.5819  &  \textbf{0.8229}  &  1.2248  &  1.1181 \\ \hline
 $\text{diag}(1,0)$& 0.9822   & 1.3890  &  \textbf{1.2248}  &  1.1181 \\ \hline
 $\text{diag}(1,1)$ & 1.1181  &  1.5812 &   \textbf{1.2248} &   1.1181 \\ \hline
\end{tabular}
\end{center}
\end{table}
\end{center}

In our next experiment, we chose different size real DH pencils $sE-(J-R)$, where 
$J=0$ is fixed, and $R$ and $E$ are randomly generated positive definite matrices. The computed distances are depicted in Table~\ref{Tab2:DHresults}. We observe that in all our experiments
(i) $\delta_0^d(J,R,E)$ gives a better upper bound than the one obtained in~\cite{MehMW20}; 
(ii) $\delta_0^d(J,R,E)$ is obtained by $\Delta_J=0$ and some rank one perturbations $\Delta_R$ and $\Delta_E$, see Remark~\ref{rem:rankonedh}. This implies that $\delta_0^d(J,R,E)={\|[\Delta_J,\Delta_R,\Delta_E]\|}_F$ and hence the last two columns in Table~\ref{Tab2:DHresults} are the same. 

\begin{center}
\begin{table}[h!]
\begin{center}
\caption{structured distances for random DH pencils $sE-(J-R)$, where $J=0$.}
\label{Tab2:DHresults}
\begin{tabular}{|c|c|c|c|c|}
\hline 
size &   $\delta_0(J,R,E)$&  $\emph{u.b.} (=\sqrt{2}\delta_0(J,R,E)$)~\cite{MehMW20}& our \emph{u.b.} &  $\delta_0^d(J,R,E)$\\ 
&\eqref{eq:unstr_DHnull} & \cite[Theorem 13]{MehMW20} &  & Theorem~\ref{thm:mainJRE} \\
\hline
$3\times 3$ &5.5970   & 7.9154    &\textbf{6.7948}  &  6.7948\\ \hline
$4\times 4$ & 1.5258   & 2.1578 &   \textbf{1.7442}  &  1.7442 \\ \hline
$5\times 5$ &  2.0962  &  2.9645  &  \textbf{2.2246}   & 2.2246\\ \hline
$6\times 6$ & 3.9478   & 5.5830  &  \textbf{4.8249} &   4.8249 \\ \hline
$7\times 7$ & 2.4811   & 3.5089&    \textbf{2.8868}  &  2.8868 \\ \hline
$8\times 8$ &  3.5417    &5.0087   & \textbf{4.6057}   & 4.6057 \\ \hline
\end{tabular}
\end{center}
\end{table}
\end{center}

}
\end{example}

\section{Conclusion}
We have proposed a purely linear algebra-based approach to derive explicit formulas to compute  structured distances to common null space for a given pencil $A+sE$ with a symmetry structure. This includes Hermitian, skew-Hermitian, $*$-even, $*$-odd, $*$-palindromic, T-palindromic, and DH pencils. For these structures, we have also obtained a family of computable lower bounds to the structured distance to singularity  for pencil $A+sE$ by deriving  a lower bound to  the nearest structured pencil that has some pre-specified $n+1$ distinct eigenvalues. 
As far as we know, this is the first work that has considered the structured distance to singularity except the one in~\cite{MehMW20} for DH pencils.

\textbf{Acknowledgment} We would like to thank Prof. Volker Mehrmann of TU Berlin, Germany,  and Prof. Shreemayee Bora of IIT Guwahati, India,  for their thoughtful suggestions for the improvement of the paper.

\small

\bibliographystyle{siam}
\bibliography{bibliostable2}

\end{document}